\documentclass[11pt]{article}
\usepackage{amsfonts}

\usepackage{CJK}
\usepackage{amsfonts,amssymb,amsmath,mathrsfs,multirow,booktabs}
\usepackage{xcolor,soul}
\usepackage{color,latexsym,amsfonts}
 \newcommand{\qed}{\hfill\rule{2mm}{3mm}\vspace{4mm}}

 \newtheorem{theorem}{Theorem}[section]
 \newtheorem{lemma}[theorem]{Lemma}
 \newtheorem{corollary}[theorem]{Corollary}
 \newtheorem{proposition}[theorem]{Proposition}
 \newtheorem{example}[theorem]{Example}
 \newtheorem{Definition}[theorem]{Definition}
 \newtheorem{remark}[theorem]{Remark}
 \newtheorem{condition}[theorem]{Condition}
 \newtheorem{conjecture}[theorem]{Conjecture}

 \def\blemma{\begin{lemma}\sl{}\def\elemma{\end{lemma}}}
 \def\btheorem{\begin{theorem}\sl{}\def\etheorem{\end{theorem}}}

 \def\bremark{\begin{remark}\sl{}\def\eremark{\end{remark}}}

 \def\beqlb{\begin{eqnarray}}\def\eeqlb{\end{eqnarray}}
 \def\beqnn{\begin{eqnarray*}}\def\eeqnn{\end{eqnarray*}}

 \def\mbf{\mathbf}

 \def\<{\langle}\def\>{\rangle}

 \def\ar{&\!\!}

 \def\eqref#1{{\rm(\ref{#1})}}

 \def\proof{\noindent{\it
 Proof.~}}\def\qed{\hfill$\Box$\medskip}

\def\e{{\mbox{\rm e}}}
\def\<{\left<}\def\>{\right>}

  \def\mbf{\mathbf}
\newcommand{\dd}{\mathrm{d}}

\newfam\msbmfam\font\tenmsbm=msbm10\textfont
\msbmfam=\tenmsbm\font\sevenmsbm=msbm7
\scriptfont\msbmfam=\sevenmsbm

\def\<{\left<}\def\>{\right>}
\def\({\left(}\def\){\right)}

\begin{document}

\centerline{\Large\bf Boundary behaviors for a class of continuous-state}
\medskip
\centerline{\Large\bf nonlinear branching  processes in critical cases}
\bigskip

\centerline{
Shaojuan Ma\footnote{School of Mathematics and Information Science,
North Minzu University, Yinchuan, China. Email: sjma@nmu.edu.cn.},
Xu Yang\footnote{School of Mathematics and Information Science,
North Minzu University, Yinchuan, China. Email: xuyang@mail.bnu.edu.cn. Corresponding author.}
and Xiaowen Zhou\footnote{Department of Mathematics and Statistics, Concordia University, Montreal, Canada.
Email: xiaowen.zhou@concordia.ca.}}

\bigskip\bigskip

{\narrower{\narrower

\noindent{\bf Abstract.}
Using Foster-Lyapunov techniques we establish new conditions on non-extinction, non-explosion, coming down from infinity and staying infinite, respectively, for the general continuous-state nonlinear branching processes introduced in Li et al. (2019). These results can be applied to
identify  boundary behaviors   for the critical cases of the above nonlinear branching processes with power rate functions  driven by Brownian motion and (or) stable Poisson random measure, which was left open in  Li et al. (2019). In particular, we show that even in the critical cases,  a phase transition happens between coming down from infinity and staying infinite.

\bigskip


\bigskip

\par}\par}

\section{Introduction}

Continuous-state branching processes (CSBPs in short) are nonnegative-valued  Markov processes satisfying the additive branching processes. They often
arise as time-population scaling limits of discrete-state branching processes, and can also be obtained from spectrally positive L\'evy processes via the Lamperti transform. The introduction of CSBP allows the applications of stochastic analysis, L\'evy processes and stochastic differential equations (SDE in short) techniques to its study.  We refer to Li (2011), Li (2019) and Kyprianou (2012) and references therein for  comprehensive reviews on CSBPs.

Generalized versions of the CSBP  have been proposed in recent years to incorporate  interactions between individuals and (or) between individuals and the population. A class of CSBPs with nonlinear branching mechanism, obtained by generalized Lamperti transform, is introduced in Li (2019). In Li et al. (2019),  a more general version of the nonlinear CSBP is proposed as the solution to SDE
\begin{equation}\label{sde}
X_t=
x+\int_0^t a_0(X_s)\dd s+\int_0^t\int_0^{a_1(X_s)}W(\dd s,\dd u)
+\int_0^t\int_{(0,\infty)}\int_0^{a_2(X_{s-})}z\tilde{M}(\dd s,\dd z,\dd u),
\end{equation}
where $x>0$,
$a_0$ and $a_1,a_2\ge0$ are Borel functions
on $[0,\infty)$, $W(\dd s,\dd u)$ and
$\tilde{M}(\dd s,\dd z,\dd u)$
denote a Gaussian white noise and an
independent compensated Poisson random measure,
respectively.
The model of Li (2019) corresponds to  solution $(X_t)_{t\geq 0}$ to SDE (\ref{sde}) with identical power rate functions $a_i, i=0,1,2$.

These nonlinear CSBPs allow  richer behaviors such as coming down from infinity.
Some extinction, extinguishing, explosion and coming down from infinity properties are proved in Li (2019).  By analyzing weighted occupation times for spectrally positive L\'evy process, asymptotic results on the speeds of coming down from infinity and explosion are obtained in Foucart et al. (2019) and Li and Zhou (2019), respectively, for nonlinear CSBP corresponding to solution to SDE (\ref{sde}) with identical rate functions $a_0=a_1=a_2$.
Exponential ergodicity  for the general continuous-state nonlinear branching processes in Li et al. (2019) is studied by Li and Wang (2020) using coupling techniques.

A version of SDE (\ref{sde}) with $a_1\equiv 0$ and power functions $a_0$ and $a_2$ is considered earlier in  Berestycki et al. (2015) where using the Lamperti  transform,  a necessary and sufficient condition for extinction is obtained and the pathwise uniqueness of solution is studied.  Work on the continuous-state logistic branching process can be found in Lambert (2005), Le et al. (2013) and Le (2014).

Using a martingale approach, the extinction, explosion and coming down from infinity behaviors are further discussed in Li et al. (2019) and some rather sharp criteria in terms of $a_0, a_1, a_2$ and $\mu$ are obtained on characterization of different kinds of boundary behaviors for the nonlinear CSBPs as a Markov process. In   Example 2.18 of Li et al. (2019) where $a_0, a_1$ and $a_2$ are taken to be power functions and $\tilde{M}$ is taken to be an $\alpha$-stable Poison random measure with index $\alpha\in (1, 2)$. The above criteria are further expressed in terms of the coefficients and the powers of functions $a_i$s and the stable index $\alpha$. But for the critical cases, where the coefficients, the powers and the index $\alpha$ satisfy certain linear equations, the martingale approach fails and the corresponding boundary classification remains an open problem.

The main goal of  this paper is to identify the exact  boundary behaviors in the above mentioned critical cases for the solution $(X)_{t\geq 0}$ to (\ref{sde}). For this purpose, we adapt the Foster-Lyapunov approach and select logarithm type test functions to obtain two new conditions under which the nonlinear CSBP never becomes extinct and never explodes, respectively. Similarly, for the boundary at infinity,  we also find a  Foster-Lyapunov condition with which we can show that an interesting phase transition occurs between coming down from infinity and staying infinite  for different choices of coefficients and powers of the power functions and differential values of the stable index $\alpha$.

The rest of the paper is arranged as follows. We introduce the generalized CSBPs with nonlinear branching in more details and present the main theorem in Section \ref{main}. The proofs of preliminary results  on Foster-Lyapunov criteria and the main theorem are deferred to Section \ref{proof}.

\section{Main results}\label{main}
\setcounter{equation}{0}

Let $U$ be a Borel set on $(0,\infty)$.
Given $\sigma$-finite measures $\mu$ and $\nu$ on $(0, \infty)$ such that
	\[(z\wedge z^2)\mu(\dd z)\,\,\text{and}\,\, (1\vee\ln(1+z))\nu(\dd z)\,\, \text{are finite measures
on $U$ and $(0,\infty)\setminus U$, respectively},\]
we consider the following SDE that is a modification of (\ref{sde}):
 \beqlb\label{1.01}
X_t
 \ar=\ar
x+\int_0^t a_0(X_s)\dd s+\int_0^t\int_0^{a_1(X_s)}W(\dd s,\dd u) \cr
 \ar\ar
+\int_0^t\int_U\int_0^{a_2(X_{s-})}z\tilde{M}(\dd s,\dd z,\dd u)
+\int_0^t\int_{(0,\infty)\setminus U}\int_0^{a_3(X_{s-})}zN(\dd s,\dd z,\dd u),
 \eeqlb
where $x>0$,
$a_0$ and $a_1,a_2,a_3\ge0$ are Borel functions
on $[0,\infty)$,
$W(\dd s,\dd u)$ is a Gaussian white noise with density $\dd s\dd u$,
$\tilde{M}(\dd s,\dd z,\dd u)$
denotes a compensated Poisson random measure on
$(0,\infty)\times U\times(0,\infty)$
with density $\dd s\mu(\dd z)\dd u$,
and
$N(\dd s,\dd z,\dd u)$
denotes a Poisson random measure on $(0,\infty)\times ((0,\infty)\setminus U)\times(0,\infty)$
with density $\dd s\nu(\dd z)\dd u$.
We assume that $W(\dd s,\dd u)$, $\tilde{M}(\dd s,\dd z,\dd u)$
and $N(\dd s,\dd z,\dd u)$ are independent of each other.
A similar SDE (\ref{1.01}) is considered in  Li et al. (2019)  under the assumption that $U=(0,\infty)$.

We only consider the solution of \eqref{1.01} before the minimum of their
first times of hitting zero and  reaching infinity
(that is $\tau_0^-$ and $\tau_\infty^+$,
which will be given in the following), respectively, i.e. both zero and  infinity are absorbing states for the solution. See Section 2 of Li et al. (2019) for more details.
By the same argument as Theorem 3.1 in Li et al. (2019), we can show that SDE \eqref{1.01} has a pathwise unique solution if functions $a_0,a_1,a_2,a_3$ are
locally Lipschitz on $(0,\infty)$.
The main purpose of this paper is to investigate the extinction, explosion and
coming down from infinity behaviors,
and the uniqueness of solution to SDE (\ref{1.01}) is not required.

Throughout this paper we always assume that
$a_0,a_1,a_2,a_3$ are bounded on any bounded interval of $[0,\infty)$\
and that process $(X_t)_{t\ge0}$ is defined on filtered probability space
$(\Omega,\mathscr{F},\mathscr{F}_t,\mbf{P})$
which satisfies the usual hypotheses.
We use $\mbf{P}_x$ to denote the law of a
process started at $x$, and denote by $\mbf{E}_x$ the associated expectation.
For $a, b>0$ we define first passage times
 \beqnn
\tau_a^-:=\inf\{t>0:X_t<a\},\quad\tau_b^+:=\inf\{t>0:X_t>b\}
 \eeqnn
and
 \beqnn
\tau_0^-:=\inf\{t>0:X_t=0\},\quad
\tau_\infty^+:=\lim_{n\to\infty}\tau_n^+
 \eeqnn
with the convention $\inf\emptyset=\infty$.
Let $C^2((0,\infty))$ denote the space of
twice continuously differentiable functions on $(0,\infty)$.

We next introduce several functions.
For $u>0$ let
\beqnn
\phi(u)
\ar:=\ar
-a_0(u)u^{-1}+\frac12a_1(u)u^{-2}
+a_2(u)\int_U z^2\mu(\dd z)\int_0^1(u+vz)^{-2}(1-v)\dd v \cr
\ar\ar
-a_3(u)\int_{(0,\infty)\setminus U} z\nu(\dd z)\int_0^1(u+vz)^{-1}\dd v.
\eeqnn
For $\rho,z>0$ and $u>3$ let
 \beqnn
K_\rho(u,z):=\Big(\frac{\ln (u+z)}{\ln u}\Big)^{-\rho}
+\rho\frac{\ln (u+z)}{\ln u}-(\rho+1)>0
 \eeqnn
and
 \beqnn
H_\rho(u):=\frac12a_1(u)u^{-2}
+a_2(u)\int_U K_\rho(u,z)\mu(\dd z)
+a_3(u)\int_{(0,\infty)\setminus U} K_\rho(u,z)\nu(\dd z).
 \eeqnn

Process $(X_t)_{t\geq 0}$ becomes extinct if $\tau^-_0<\infty$; it explodes if $\tau^+_\infty<\infty$; it stays infinite if
$\lim_{x\to\infty}\mbf{P}_x\{\tau_a^-<t\}=0$ for all $t>0$ and all large $a$; it comes down from infinity if
$\lim_{a\to\infty}\lim_{x\to\infty}\mbf{P}_x\{\tau_a^-<t\}=1$ for all $t>0$.

The following main result provides new criteria on non-extinction, non-explosion, coming down from infinity and staying infinite for the solution to SDE (\ref{sde}).
 \btheorem\label{t5.6} For the solution $(X_t)_{t\geq 0}$ to (\ref{1.01}) we have
 \begin{itemize}
\item[{(i)}] if
$\phi(u)\le 0$ for all small enough $u>0$,
then $\mbf{P}_x\{\tau_0^-=\infty\}=1$ for all $x>0$, i.e. there is no extinction;
\item[{(ii)}] if
$\phi(u)\ge 0$ for all large enough $u>0$,
then $\mbf{P}_x\{\tau_\infty^+=\infty\}=1$ for all $x>0$, i.e. there is no explosion;
\item[{(iii)}]  if
$\phi(u)\le 0$ for all large enough $u>0$ and
\[\limsup_{u\rightarrow \infty}H_\rho(u)<\infty\]
for some constant $\rho>0$,
then the process $(X_t)_{t\ge0}$ stays infinite;
\item[{(iv)}] if $\phi(u)\ge 0$ for all large enough $u>0$ and
 \beqnn
\lim_{u\rightarrow\infty}(\ln u)^{-\rho-2} H_\rho(u)=\infty
 \eeqnn
for some constant $\rho>0$,
then the  process $(X_t)_{t\ge0}$ does not explode and it comes down from infinity.
\end{itemize}
 \etheorem

\begin{remark}
Since the process $(X_t)_{t\ge0}$ does not explode and comes down from infinity	under the assumptions of Theorem \ref{t5.6} (iv),  then under additional assumptions on functions $a_i$s, $(X_t)_{t\ge0}$ can be extended to a Feller process defined on state space $[0, \infty]$ with $\infty$ as  its entrance boundary; see Foucart et al. (2020).


\end{remark}	
	
Until the end of this section we focus on the special  case that  $U=(0,\infty)$, $a_0,a_1,a_2$ are power functions and $\mu(\dd z)$ is an $\alpha$-stable measure,
that is
\begin{equation}\label{con_A}
a_i(u)=b_i u^{r_i}~ \text{for}~ i=0, 1, 2 ~\text{with}~r_0, r_1, r_2\geq 0, b_1,b_2\ge0, b_0>0
\end{equation}
and
 \begin{equation}\label{con_B}
\mu(\dd z)=\frac{\alpha(\alpha-1)}{\Gamma(2-\alpha)}1_{\{z>0\}}z^{-1-\alpha}\dd z~\text{for Gamma function}\, \Gamma \,\text{and}\, 1<\alpha<2.
\end{equation}

By the properties of Gamma function we have
 \beqnn
\int_0^\infty z^2\mu(\dd z)\int_0^1(u+vz)^{-2}(1-v)\dd v=\Gamma(\alpha)u^{-\alpha}
 \eeqnn
and then
 \beqnn
\phi(u)=-b_0u^{r_0-1}+\frac12b_1u^{r_1-2}+\Gamma(\alpha)b_2u^{r_2-\alpha}, \quad u>0.
 \eeqnn

We further estimate
\[H_\rho(u)=\frac12a_1(u)u^{-2}+a_2(u)\int_0^\infty K_\rho(u,z)\mu(\dd z)\]
for which we first estimate $\int_0^\infty K_\rho(u,z)\mu(\dd z)$.
Note that for $y>0$ and $f(y):=y^{-\rho}+\rho y-(\rho+1)$, by Taylor's formula we have
 \beqnn
f(y)
 \ar=\ar
f(1+y-1)=f(1+y-1)-f(1)-(y-1)f'(1) \cr
 \ar=\ar
(y-1)^2\int_0^1f''(1+v(y-1))(1-v)\dd v.
 \eeqnn
Then by a change of variable,
 \beqnn
 \ar\ar
\int_0^\infty K_\rho(u,z)\mu(\dd z) \cr
 \ar=\ar
\rho(\rho+1)\int_0^\infty\mu(\dd z)
\int_0^1\Big(\frac{\ln (u+vz)}{\ln u}-1\Big)^2\Big(1+\frac{v\ln(u+vz)}{\ln u}-v\Big)^{-\rho-2}(1-v)\dd v \cr
 \ar=\ar
\rho(\rho+1)u^{-\alpha}\int_0^\infty\mu(\dd z)
\int_0^1\Big(\frac{\ln (1+vz)}{\ln u}\Big)^2\Big(1+\frac{v\ln(1+vz)}{\ln u}\Big)^{-\rho-2}(1-v)\dd v.
 \eeqnn
Observe that
$\ln (1+z)\le C(z\wedge \sqrt{z})$
for all $z>0$ and some constant $C>0$,
which implies that for $u>3$,
 \beqlb\label{1.2}
\int_0^\infty K_\rho(u,z)\mu(\dd z)
 \ar\le\ar
\rho(\rho+1)u^{-\alpha} (\ln u)^{-2}
\int_0^\infty(\ln (1+z))^2\mu(\dd z) \cr
 \ar\le\ar
C^2\rho(\rho+1)u^{-\alpha} (\ln u)^{-2}\int_0^\infty (z\wedge z^2)\mu(\dd z).
 \eeqlb
Moreover, it is elementary to see that
 \beqlb\label{1.3}
\ar\ar\int_0^\infty K_\rho(u,z)\mu(\dd z)\cr
\ar\ge\ar
\rho(\rho+1)u^{-\alpha}(\ln u)^{-2}\Big(1+\frac{\ln3}{\ln u}\Big)^{-\rho-2}\int_1^2 (\ln (3/2))^2\mu(\dd z)
\int_{1/2}^1(1-v)\dd v.
 \eeqlb


	\begin{remark}
		In Section 2.5 of Li et al. (2019),
		the exact conditions are found for the above mentioned model with polynomial rate functions to exhibit  extinction/non-extinction, explosion/non-explosion and coming-down-from-infinity/staying-infinite behaviors, respectively, except for the critical case that
\[b_0 =\frac{b_1}{2} +\Gamma(\alpha) b_2>0, \, r_1=r_0+1\,\,\text{ when}\,\, b_1>0\,\, \text{ and} \,\,
r_2=r_0+\alpha-1 \,\, \text{ when} \,\, b_2>0.\]
Observe that in this critical case $\phi(u)=L(\ln)(u)=0$, where operator $L$, to be defined in (\ref{2.1}), denotes the generator of process $X$. This inspires us to choose  logarithm type test functions for the main proofs.
	\end{remark}

As the main goal of this paper, applying Theorem \ref{t5.6} together with \eqref{1.2}--\eqref{1.3},  we provide an answer to this open problem.
\begin{corollary}\label{example}
Suppose that (\ref{con_A}) and (\ref{con_B}) hold with $b_0 =\frac{b_1}{2} +\Gamma(\alpha) b_2>0$,
 $r_1=r_0+1$ when $b_1>0$ and
$r_2=r_0+\alpha-1$ when $b_2>0$.
Then $\phi(u)=0$ for all $u>0$, and we have for all $x>0$,
 \beqnn
\mbf{P}_x\{\tau_0^-=\infty\}=1\quad\text{and}\quad
\mbf{P}_x\{\tau_\infty^+=\infty\}=1.
 \eeqnn
Moreover,
process $(X_t)_{t\ge0}$ stays infinite if both $r_1\le 2$ (when $b_1>0$) and $r_2\le \alpha$ (when $b_2>0$),
and it comes down from infinity
if either $r_1> 2$ (when $b_1>0$) or $r_2> \alpha$ (when $b_2>0$).
\end{corollary}

\bremark
Note that in the critical cases, there is an interesting phase transition between coming down from infinity and staying infinite. Intuitively, in these cases the process comes down from infinity if the fluctuations are relatively large and stays infinite otherwise.

Combining Corollary \ref{example} and Example 2.18 of Subsection 2.5 in Li et al. (2019), we recover the necessary and sufficient condition on the extinction of   solution to  the SDE of
 Berestycki et al. (2015); see Theorem 1.1 there.
\eremark

\section{Proofs}\label{proof}
\setcounter{equation}{0}

Before presenting the proof of Theorem \ref{t5.6} we first prove some preliminary Foster-Lyapunov criteria.
Suppose that $g\in C^2((0,\infty))$ satisfies
 \beqlb\label{2.5}
\sup_{z\ge1,u\ge v}\big[|g'(u)|+|g''(u)|+ |g(u+z)-g(u)|/\ln(1+z)\big]<\infty
 \eeqlb
for all $v>0$.
For $u>0$, put
 \beqlb\label{2.1}
Lg(u)
 \ar:=\ar
a_0(u)g'(u)+\frac12a_1(u)g''(u)
+a_2(u)\int_{U}[g(u+z)-g(u)-zg'(u)]\mu(\dd z) \cr
 \ar\ar
+a_3(u)\int_{(0,\infty)\setminus U} [g(u+z)-g(u)]\nu(\dd z)\cr
 \ar=\ar
a_0(u)g'(u)+\frac12a_1(u)g''(u)
+a_2(u)\int_U z^2\mu(\dd z) \int_0^1g''(u+zv)(1-v)\dd v \cr
 \ar\ar
+a_3(u)\int_{(0,\infty)\setminus U}z\nu(\dd z)\int_0^1g'(u+zv)\dd v
 \eeqlb
by Taylor's formula.
By It\^o's formula,
 \beqnn
g(X_t)
 \ar=\ar
g(x)+\int_0^tLg(X_s)\dd s
+\int_0^t\int_U\int_0^{a_2(X_{s-})}[g(X_{s-}+z)-g(X_{s-})]\tilde{M}(\dd s,\dd z,\dd u) \cr
 \ar\ar
+\int_0^t\int_{(0,\infty)\setminus U}\int_0^{a_3(X_{s-})}
[g(X_{s-}+z)-g(X_{s-})]\tilde{N}(\dd s,\dd z,\dd u)
 \eeqnn
For $b>a>0$ let $\gamma_{a,b}:=\tau_{a}^-\wedge\tau_b^+$
and
 \beqnn
M_t^g:=
g(X_t)-g(x)
-\int_0^tLg(X_s)\dd s.
 \eeqnn
Then under condition \eqref{2.5},
 \beqlb\label{2.2}
t\mapsto M_{t\wedge \gamma_{a,b}}^g
\quad \mbox{is a martingale}
 \eeqlb
for all $b>a>0$.

\blemma\label{t2.2}
Given $0<a<x<b<\infty$, for any function $g\in C^2((a, b))$
satisfying
\eqref{2.5} and constant $d_{a,b}>0$ satisfying
 \beqnn
Lg(u)\le d_{a,b}g(u),\qquad u\in(a,b),
 \eeqnn
we have
 \beqlb\label{2.3}
\mbf{E}_x\big[g(X_{t\wedge \gamma_{a,b}})\big]
\le
g(x)\e^{d_{a,b}t},\qquad t\geq 0.
 \eeqlb
\elemma
\proof
It follows from \eqref{2.2} that
 \beqnn
\mbf{E}_x\big[g(X_{t\wedge \gamma_{a,b}})\big]=g(x)
+\mbf{E}_x\Big[\int_0^{t\wedge \gamma_{a,b}}Lg(X_s)\dd s\Big]
\le
g(x)+d_{a,b}\int_0^t\mbf{E}_x\big[g(X_{s\wedge \gamma_{a,b}})\big]\dd s.
 \eeqnn
By Gronwall's lemma,
 \beqnn
\mbf{E}_x\big[g(X_{t\wedge \gamma_{a,b}})\big]
\le
g(x)\e^{d_{a,b}t},
 \eeqnn
which ends the proof.
\qed

\blemma\label{t2.1}
Let $(X_t)_{t\geq 0}$ be the solution to SDE (\ref{1.01}).
\begin{itemize}
\item[{\normalfont(i)}]
For any fixed $b>0$,
if there exists a function $g\in C^2((0, \infty))$ strictly positive on $(0,b]$  satisfying \eqref{2.5} and $\lim_{u\to0}g(u)=\infty$,
and there is a constant $d(b)>0$ such that
$Lg(u)\le d(b)g(u)$ for all $0<u<b$,
then $\mbf{P}_x\{\tau_0^-\ge\tau^+_b\}=1$ for all $0<x<b$.

\item[{\normalfont(ii)}]
For any fixed $a>0$,
if there exists a function $g\in C^2((0,\infty))$ strictly positive on $[a,\infty)$ satisfying
\eqref{2.5} and $\lim_{u\to\infty}g(u)=\infty$,
and there is a constant $d(a)>0$ such  that
$Lg(u)\le d(a)g(u)$ for all $u>a$,
then $\mbf{P}_x\{\tau_\infty^+\ge\tau_a^-\}=1$ for all $x>a$.

\item[{\normalfont(iii)}]
If there exists a function $g\in C^2((0,\infty))$ strictly positive on $[u,\infty)$ for all large $u$ satisfying
\eqref{2.5} and $\lim_{u\to\infty}g(u)=0$,
 and for any large $a>0$, there is a constant $d(a)>0$ such that
$Lg(u)\le d(a)g(u)$ for all $u>a$,
then $(X_t)_{t\ge0}$ stays infinite.
\end{itemize}
\elemma
\proof
We apply Lemma \ref{t2.2} for the proofs.

For part (i), \eqref{2.3} holds for all $0<a<b$ and with
$d_{a,b}$ replaced by $d(b)$.
Then using Fatou's lemma we have
 \beqnn
\mbf{E}_x\big[\liminf_{a\to0}g(X_{t\wedge \tau_a^-\wedge\tau_b^+})\big]
\le
\liminf_{a\to0}\mbf{E}_x\big[g(X_{t\wedge \gamma_{a,b}})\big]
\le g(x)\e^{d(b)t}.
 \eeqnn
Since $\lim_{u\to0}g(u)=\infty$,
then
$\mbf{P}_x\{\tau_0^->t\wedge\tau_b^+\}=1$ for all $t,b>0$.
Letting $t\to\infty$ we obtain $\mbf{P}_x\{\tau_0^-\ge\tau_b^+\}=1$, which gives the first assertion.

For part (ii), \eqref{2.3} holds for all $b>a$ and
with $d_{a,b}$ replaced by $d(a)$.
Then using Fatou's lemma again we obtain
 \beqnn
\mbf{E}_x\big[\liminf_{b\to\infty}g(X_{t\wedge \tau_a^-\wedge\tau_b^+})\big]
\le
\liminf_{b\to\infty}\mbf{E}_x\big[g(X_{t\wedge \gamma_{a,b}})\big]
\le g(x)\e^{d(a)t}.
 \eeqnn
Since $\lim_{u\to\infty}g(u)=\infty$, then
$\mbf{P}_x\{\tau_\infty^+>t\wedge\tau_a^-\}=1$ for all $t>0$.
Letting $t\to\infty$ we get
$\mbf{P}_x\{\tau_\infty^+\ge\tau_a^-\}=1$,
which implies the second assertion.

For part (iii), given any large $a>0$, \eqref{2.3} holds for all $b>a$ and
with $d_{a,b}$ replaced by $d(a)$ again.
We can also get
 \beqnn
\mbf{E}_x\big[g(X_{\tau^-_a})1_{\{\tau^-_a<t\wedge\tau_\infty^+\}}\big]
\ar\le\ar
\liminf_{b\to\infty}
\mbf{E}_x\big[g(X_{\tau^-_a})1_{\{\tau^-_a<t\wedge\tau_b^+\}}\big]\cr
\ar\le\ar
\liminf_{b\to\infty}\mbf{E}_x\big[g(X_{t\wedge \gamma_{a,b}})\big]
\le g(x)\e^{d(a)t},
 \eeqnn
which implies
 \beqnn
g(a)\mbf{P}_x\{\tau^-_a<t\wedge\tau_\infty^+\}
\le g(x)\e^{d(a)t}.
 \eeqnn
Since $\lim_{u\to\infty}g(u)=0$, then for all $t,a>0$,
 \beqlb\label{2.4}
\lim_{x\to\infty}\mbf{P}_x\{\tau^-_a<t\wedge\tau_\infty^+\}=0.
 \eeqlb
Observe that
$\{\tau^-_a\ge\tau_\infty^+\}\subset\{\tau^-_a=\infty\}$.
Then combining \eqref{2.4} we have
 \beqnn
\lim_{x\to\infty}\mbf{P}_x\{\tau^-_a<t\}
\le
\lim_{x\to\infty}\mbf{P}_x\{\tau^-_a<t\wedge\tau_\infty^+\}
+\lim_{x\to\infty}\mbf{P}_x\{\tau^-_a<t,\tau_\infty^+\le\tau^-_a\}=0
 \eeqnn
for all $t>0$. Then the process stays infinite.
\qed

The next lemma provides a condition that associates  the probability of coming down from infinity with the probability of non-explosion. Its proof is a modification of Proposition 2.2 in Ren et al. (2019).

\blemma\label{t2.3}
Suppose that there exist a function $g(u)\in C^2((0,\infty))$ bounded and strictly positive for all large $u$, satisfying
\eqref{2.5} and $\limsup_{u\to\infty} g(u)>0$,
and a strictly positive function $d$ on $(0, \infty)$
such  that
\[Lg(u)\ge d(a) g(u)\,\,\text{ for all large}\,\, u\,\,\text{ and}\,\, \lim_{a\to\infty}d(a)=\infty.\]
Then for any $t>0$
 \beqnn
\lim_{a\to\infty}\lim_{x\to\infty}\mbf{P}_x\{\tau_a^-<t\}\ge
\liminf_{x\to\infty}\mbf{P}_x\{\tau_\infty^+=\infty\}.
 \eeqnn
Consequently, process $(X_t)_{t\ge0}$ comes down from infinity if there is no explosion.
\elemma
\proof
The proof is a modification of that of Proposition 2.2 in Ren et al. (2019).
We present the details for completeness.
By \eqref{2.2}, for all large $a<b$
 \beqnn
\mbf{E}_x\big[g(X_{t\wedge\gamma_{a,b}})\big]
=
g(x)
+\mbf{E}_x\Big[\int_0^{t\wedge\gamma_{a,b}}Lg(X_s)\dd s\Big]
=
g(x)
+\int_0^t \mbf{E}_x\Big[Lg(X_s)1_{\{s\le \gamma_{a,b}\}}\Big] \dd s
 \eeqnn
and then by integration by parts,
 \beqnn
 \ar\ar
\mbf{E}_x\big[g(X_{t\wedge\gamma_{a,b}})\big]\e^{-d(a)t} \cr
 \ar=\ar
g(x)+\int_0^t \mbf{E}_x\big[g(X_{s\wedge\gamma_{a,b}})\big]\dd (\e^{-d(a)s})
+\int_0^t \e^{-d(a)s}\dd \big(\mbf{E}_x\big[g(X_{s\wedge\gamma_{a,b}})\big]\big) \cr
\ar=\ar
g(x)
-d(a)\int_0^t\mbf{E}_x\big[g(X_{s\wedge\gamma_{a,b}})\e^{-d(a)s}\big]\dd s
+\int_0^t\e^{-d(a)s}\mbf{E}_x\big[Lg(X_s)1_{\{s\le\gamma_{a,b}\}}\big]\dd s \cr
\ar\ge\ar
g(x)
-d(a)\int_0^t\mbf{E}_x\big[g(X_{s\wedge\gamma_{a,b}})\big]\e^{-d(a)s}\dd s
+d(a)\int_0^t\e^{-d(a)s}\mbf{E}_x\big[g(X_s)1_{\{s\le\gamma_{a,b}\}}\big]\dd s,
 \eeqnn
which implies that
 \beqnn
g(x)\le\mbf{E}_x\big[g(X_{t\wedge\gamma_{a,b}})\e^{-d(a)t}\big]
+d(a)\mbf{E}_x\Big[\int_0^tg(X_{\gamma_{a,b}})\e^{-d(a)s}
1_{\{s>\gamma_{a,b}\}}\dd s\Big].
\eeqnn
Letting $t\to\infty$ in the above inequality and using the dominated convergence we obtain
 \beqnn
g(x)
\le
d(a)\mbf{E}_x\Big[g(X_{\gamma_{a,b}})\int_{\gamma_{a,b}}^\infty\e^{-d(a)s}\dd s\Big]
=\mbf{E}_x\big[g(X_{\gamma_{a,b}}) \e^{-d(a)\gamma_{a,b}}\big].
 \eeqnn
It follows that
 \beqnn
g(x)
 \ar\le\ar
\mbf{E}_x\Big[\lim_{b\to\infty}g(X_{\gamma_{a,b}})
\e^{-(\tau_a^-\wedge\tau^+_\infty)d(a)}
\big(1_{\{\tau_\infty^+<\tau^-_a\}}+1_{\{\tau_a^-<t,\tau_a^-\le\tau^+_\infty\}}
+1_{\{t\le\tau^-_a\le\tau_\infty^+\}}\big)\Big] \cr
 \ar\le\ar
\limsup_{u\to\infty}g(u)
\mbf{P}_x\{\tau_\infty^+<\infty\}
+g(a)\mbf{P}_x\{\tau_a^-<t,\tau_a^-\le\tau^+_\infty \}+g(a)\e^{-d(a) t} \cr
 \ar\le\ar
\limsup_{u\to\infty}g(u)( 1-\mbf{P}_x\{\tau_\infty^+=\infty\})
+g(a)\mbf{P}_x\{\tau_a^-<t\}+g(a)\e^{-d(a) t}.
 \eeqnn
Letting $x\to\infty$ first,
 \beqnn
\limsup_{x\to\infty}g(x)
 \ar\leq\ar
\limsup_{u\to\infty}g(u)
\limsup_{x\to\infty}(1-\mbf{P}_x\{\tau_\infty^+=\infty\}) \cr
 \ar\ar
+g(a)\lim_{x\to\infty}\mbf{P}_x\{\tau_a^-<t \}+g(a)\e^{-d(a) t}.
 \eeqnn
Then letting $a\to\infty$, by the conditions in the lemma we have
 \beqnn
\limsup_{x\to\infty}g(x)
 \ar\leq\ar
\limsup_{u\to\infty}g(u)
\left(1-\liminf_{x\to\infty}\mbf{P}_x\{\tau_\infty^+=\infty\}\right)\cr
 \ar\ar
+\limsup_{a\to\infty}g(a)\limsup_{a\to\infty}
\lim_{x\to\infty}\mbf{P}_x\{\tau_a^-<t \}.
 \eeqnn
Observing that $\lim_{x\to\infty}\mbf{P}_x\{\tau_a^-<t \} $ is increasing in $a$, the desired inequality then follows from the above inequality.
\qed


We are now ready to show the proofs of the main results.

\noindent{\it Proof of Theorem \ref{t5.6}}.
(i) Suppose that there is a constant $0<c_1<1$ so that
$\phi(u)\le 0$ for all $0<u<c_1$.
For $n\ge1$ let $g_n(u)=1+\ln n+\ln u^{-1}$.
Then $g_n(u)>0$ for $0<u\le n$
and $Lg_n(u)=\phi(u)$ by \eqref{2.1}.
Thus $Lg_n(u)\le 0$ for $0<u<c_1$.
Since $a_0,a_1,a_2,a_3$ are bounded on $[c_1,n]$,
$Lg_n$ is bounded on $[c_1,n]$.
Now using Lemma \ref{t2.1}(i) we obtain $\mbf{P}_x\{\tau_0^-\ge\tau^+_n\}=1$
for all $0<x<n$.
Since the process is defined before the first time
of hitting zero or explosion,
$\mbf{P}_x\{\tau_0^-=\infty$ or $\tau_\infty^+=\infty\}=1$. Letting $n\to\infty$ we prove the assertion.

(ii) Suppose that there is a constant $c_2>1$ so that
$\phi(u)\ge 0$ for all $u> c_2$. Let $g_n(u)=\ln u+\ln n+1$ for $n\ge1$.
Then $g_n(u)\ge1$ for $u\ge n^{-1}$.
It follows from \eqref{2.1} that
$Lg_n(u)=-\phi(u)$ for all $u\ge n^{-1}$.
Then for all $n\ge1$,
$Lg_n(u)\le0$ for all $u\ge c_2$
and $Lg_n$ is bounded on $[n^{-1},c_2]$,
which gives $\mbf{P}_x\{\tau_\infty^+>\tau_{1/n}^-\}=1$ for all $x>n^{-1}$ by Lemma \ref{t2.1}(ii).
Letting $n\to\infty$ we have $\mbf{P}_x\{\tau_\infty^+>\tau_0^-\}=1$ for all $x>0$.
The assertion for  (ii) then follows from the definition of the solution to SDE (\ref{1.01}).

(iii)
Suppose that there exist constants $c_3>3$ and $c_4>0$ so that
$\phi(u)\le 0$ and $H_\rho(u)\le c_4$ for all $u> c_3$.
Let $g\in C^2((0,\infty))$ be a strictly positive function
with $g(u)=(\ln u)^{-\rho}$ for $\rho>0$
and $u>3$.

Then for $u>3$,
  \beqnn
 g(u+z)-g(u)
 \ar=\ar
 -\rho (\ln u)^{-\rho-1}[\ln(u+z)-\ln u]+(\ln u)^{-\rho}K_\rho(u,z) \cr
 \ar=\ar
 -\rho (\ln u)^{-\rho-1}z\int_0^1(u+vz)^{-1}\dd v+(\ln u)^{-\rho}K_\rho(u,z)
 \eeqnn
 and
 \beqnn
g'(u)=-\rho (\ln u)^{-\rho-1}u^{-1},~
g''(u)=\rho (\ln u)^{-\rho-1}u^{-2}+\rho(\rho+1) (\ln u)^{-\rho-2}u^{-2}.
 \eeqnn
Consequently, for all $u>3$ and $z>0$,
 \beqnn
g(u+z)-g(u)-zg'(u)
 \ar=\ar
-\rho(\ln u)^{-\rho-1}\big[\ln (u+z)-\ln u-zu^{-1}\big]
+(\ln u)^{-\rho}K_\rho(u,z) \cr
 \ar=\ar
\rho (\ln u)^{-\rho-1}z^2\int_0^1(u+zv)^{-2}(1-v)\dd v
+(\ln u)^{-\rho}K_\rho(u,z).
 \eeqnn
It follows that
 \beqnn
Lg(u)
 \ar=\ar
\rho (\ln u)^{-\rho-1}\phi(u)
+\frac12\rho(\rho+1) (\ln u)^{-2}a_1(u)u^{-2}g(u) \cr
 \ar\ar
+g(u)a_2(u)\int_{U} K_\rho(u,z)\mu(\dd z)
+g(u)a_3(u)\int_{(0,\infty)\setminus U} K_\rho(u,z)\nu(\dd z) \cr
 \ar\le\ar
\rho (\ln u)^{-\rho-1}\phi(u)
+[\rho(\rho+1)+1]g(u) H_\rho(u),\qquad u>3.
 \eeqnn
Then $Lg(u)\le c_4 [\rho(\rho+1)+1] g(u)$ for all  $u>c_3$.
Thus, $(X_t)_{t\ge0}$ stays infinite for all $x>0$ by Lemma \ref{t2.1}(iii).

(iv)
Let $g\in C^2((0,\infty))$ be a bounded and strictly positive function
with $g(u)=1+(\ln u)^{-\rho}$ for $\rho>0$ and $u>3$.
It follows from the argument in (iii) that for $u>3$,
 \beqnn
Lg(u)
 \ar=\ar
\rho (\ln u)^{-\rho-1}\phi(u)
+\frac12\rho(\rho+1) (\ln u)^{-\rho-2}a_1(u)u^{-2} \cr
 \ar\ar
+(\ln u)^{-\rho}a_2(u)\int_{U} K_\rho(u,z)\mu(\dd z)
+ (\ln u)^{-\rho}a_3(u)\int_{(0,\infty)\setminus U} K_\rho(u,z)\nu(\dd z) \cr
 \ar\ge\ar
\rho (\ln u)^{-\rho-1}\phi(u)+(\rho\wedge1)(\ln u)^{-\rho-2} H_\rho(u).
 \eeqnn
Then we can conclude the proof by the assumptions for this part together with Theorem \ref{t5.6}(ii) and Lemma \ref{t2.3}.
\qed

{\bf Acknowledgements.}
This work was supported by
NSFC (Nos.~11772002 and 11771018), Major research project for North Minzu University (No. ZDZX201902) and NSERC (RGPIN-2016-06704).

\end{document}